\documentclass[acmtoms]{notepolo}

\usepackage{graphicx}
\usepackage{textpos}
\usepackage{xcolor}
\definecolor{betterBrown}{rgb}{0.6667,0.2510,0}

\newcommand{\cF}{{\mathcal F}}
\newcommand{\e}{{\mathrm{e}}}
\acmVolume{}
\acmNumber{99}
\acmYear{06}
\acmMonth{}

\usepackage{multirow}
\title{Spectral algorithms for reaction-diffusion equations}

\author{R. V. CRASTER\\ Imperial College London\\
\and R. SASSI\\ Universit\`{a} di Milano, Italy}
\begin{abstract}

A collection of codes (in MATLAB \& Fortran 77), and examples, for
solving reaction-diffusion equations in one and two space
dimensions is presented. In areas of the mathematical community
spectral methods are used to remove the stiffness associated with
the diffusive terms in a reaction-diffusion model allowing
explicit high order timestepping to be used. This is particularly
valuable for two (and higher) space dimension problems. Our aim
here is to provide codes, together with examples, to allow
practioners to easily utilize, understand and implement these
ideas; we incorporate recent theoretical advances such as
exponential time differencing methods and provide timings and
error comparisons with other more standard approaches.

The examples are chosen from the literature to illustrate points
and queries that naturally arise.

\end{abstract}

\begin{document}



\category{G. 1.8 }{Numerical Analysis}{Partial Differential
Equations-Spectral Methods} \terms{Algorithms }


\setcounter{page}{1}

\begin{bottomstuff}
This work was supported in part by the EPSRC through Grant number
GR/S47663/01.

Authors' address: R. V. Craster, Department of Mathematics,
Imperial College of Science, Technology and Medicine, London, SW7
2BZ, U.K.; Roberto Sassi, Dipartimento di Tecnologie
dell'Informazione, Universit\`{a} di Milano, via Bramante 65,
26013, Crema (CR), Italy (e-mail: sassi@dti.unimi.it)
\newline
\end{bottomstuff}

\maketitle

\begin{textblock*}{\textwidth}(-0.4cm,-12.7cm)%
\textcolor{betterBrown}{\small{\sffamily The codes can be downloaded from:} {\texttt{https://github.com/roesassi/SpectralCodes}}}%
\end{textblock*}
\vspace{-3.5mm}

\section{Introduction}
The aim of this paper is to provide a suite of practically useful
and versatile spectral algorithms (in both Matlab and Fortran 77) to
efficiently solve, numerically, systems of partial differential
equations of the general form:
\begin{equation}
 u_t=u_{xx}+u_{yy}+ f(u,v),\qquad v_t=\epsilon (v_{xx}+ v_{yy})+ g(u,v)
\label{eq:govern}
\end{equation}
where $g(u,v)$ and $f(u,v)$ are nonlinear functions of $u,v$ and
$\epsilon$ a constant; we could also consider them to additionally
be functions of the first derivatives of $u,v$, but we shall not
do so here. The methods we describe are applicable to higher
dimensions and further coupled equations, however we shall
restrain ourselves to consider two space dimensions and two
coupled equations. Such equations abound in mathematical biology,
ecology, physics and chemistry, and many wonderful mathematical
patterns and phenomena exist for special cases. In one space
dimension: non-dimensional models of epidemiology and mathematical
biology, such as,
\begin{equation}
u_t=u_{xx}+u(v-\lambda), \quad v_t=\epsilon v_{xx}-uv
\label{eqn:epidemic}
\end{equation}
emerge from the modelling. In this example, there are parameters
$\lambda, \epsilon$, where $u,v$ are the infectives and
susceptibles respectively. The parameters $\lambda$, $\epsilon$
are a removal rate and a diffusivity. These types of equation
typically yield travelling waves, \cite{murray93a}, similar in
theory and spirit to that for Fisher's equation which is the usual
paradigm. Here issues such as the speed selection for travelling
waves,  and also accelerating travelling waves are of interest;
these are relevant when one species consumes or invades another;
\cite{fisher37a} was originally interested in the propagation of
advantageous genes. Reaction diffusion equations also lead to many
other interesting phenomena, such as, pulse splitting and
shedding, reactions and competitions in excitable systems, and
stability issues. Later, we shall explicitly illustrate the
versatility of the scheme presented here versus a wide range of
examples from the literature.

Our primary interest is not really in 1D simulations, these are
relatively easily undertaken using either method of lines coupled
with spatial adaptive schemes \cite{blom94a}, or finite element
collocation schemes \cite{688}, or even just simple
Crank-Nicholson and finite difference schemes \cite{sherratt97a}.
Of these, the adaptive schemes seem preferable, in general, since
they cluster the grid points in areas of sharp solution gradients.
As such the Blom \& Zegeling code has been much utilized by one of
the authors in a wide variety of different areas
(\cite{balmforth99a,balmforth00a,craster00a}). The Matlab spectral
code we develop, for one dimension, is given in
Appendix~\ref{app:code}, and is certainly competitive with all
these schemes and often faster and easier to use.
\nocite{blom94a,sherratt97a}

Unfortunately, in two space dimensions, simulations based upon the
more conventional ideas become more time consuming
\cite{pearson93a,muratov01a}, the latter requiring an hour or so
of runtime on an SGI-Cray parallel supercomputer, although the
simulations are certainly possible and accurate. It may be that
because of this expensive simulation time, comparisons with
two-dimensional simulations appear less prevalent in the
literature than the 1D cases, even though they are certainly
important in modelling pattern creation.
\nocite{pearson93a,muratov01a} However, an idea well-known in the
spectral methods community can be used to remove the stiffness
often associated with reaction diffusion equations, and thereby
allow much larger timesteps to be utilized with an explicit
timesolver; consequently the codes run very quickly even on a
standard PC or laptop. We utilize Fast Fourier Transforms in space
utilizing an integrating factor to remove the stiff terms. The
time stepping could be just a standard explicit Runge--Kutta
method, and we initially utilize this method. Later we discuss
some refinements that adjust standard Runge-Kutta to take into
account modifications motivated by the spectral scheme.

Other recent articles on numerically simulating two dimensional
reaction-diffusion models utilize wavelet \cite{cai98a}) or high
order finite difference schemes \cite{liao02a}; the spectral
method could be thought of as the logical extension of finite
differences to infinite order, a viewpoint advanced by
\cite{fornberg98a}. Alternatively physicists and mathematicians
interested in the actual processes involved, or underlying
mathematical phenomena, use operator splitting methods
\cite{ramos02a} (see also Appendix~\ref{app:adi}) or often retreat
to finite element simulations (\cite{tang93a}) or PDETWO
\cite{melgaard81a,davidson97a} another collocation based scheme.
There is only one article that we have found \cite{jones96a} that
promotes the spectral viewpoint for reaction diffusion equations,
and we agree wholeheartedly with its philosophy.
\nocite{cai98a,liao02a}

Emboldened by the recent article of \cite{weideman01a}, and the
books by \cite{trefethen01a}, \cite{boyd01a}, and our experience
with using Matlab for other two-dimensional PDE simulations in
other contexts (\cite{balmforth04a}), we primarily utilize Matlab
as our numerical vehicle, although for comparative purposes we
also coded the routines in Fortran 77 using Fast Fourier Transform
routines from \cite{fornberg98a}, \cite{canuto88a}. Taking
advantage of the built-in routines in Matlab, the resulting Matlab
codes are extremely concise, typically a page long; an example is
given in Appendix~\ref{app:code} and several are in the
accompanying electronic files. Matlab routines are also highly
portable between different platforms; our routines require Matlab
5 or higher. \nocite{fornberg98a,boyd01a,trefethen01a,weideman01a}

\section{Formulation}
Spectral methods are extremely valuable for generating numerical
methods in almost all areas of mathematics. The ability to
generate spectrally accurate spatial derivatives means that there
is simply no excuse to differentiate poorly. When this is coupled
with Fast Fourier transforms and an elegant high-level language
such as Matlab it becomes possible to generate versatile and
powerful codes. The article \cite{weideman01a} and book
\cite{trefethen01a} are quite inspirational in showing the range
of what is possible with this combination of tools. We shall
present the theory in two space dimensions:

The integrating factor approach that we utilize is to Fourier
transform equations (\ref{eq:govern}) to obtain
\begin{equation}
 U_t(\omega_x,\omega_y,t)=
 -(\omega_x^2+\omega_y^2)U(\omega_x,\omega_y,t)
 +\cF[f(u(x,y,t),v(x,y,t))],
\label{eq:original_odeU}
\end{equation}
\begin{equation}
 V_t(\omega_x,\omega_y,t)=
 -\epsilon(\omega_x^2+\omega_y^2)V(\omega_x,\omega_y,t)
 +\cF[g(u(x,y,t),v(x,y,t))]
\label{eq:original_odeV}
\end{equation}
where $U, V$ are the double Fourier transforms of $u,v$, that is,
\begin{equation}
 \cF[u(x,y,t)]= U(\omega_x,\omega_y,t)=
 \int_{-\infty}^\infty \int_{-\infty}^\infty
 u(x,y,t)\e^{-i(\omega_x x+\omega_y y)} dx dy.
\end{equation}
Let us set $\Omega^2=\omega_x^2+\omega_y^2$, and explicitly remove
the linear pieces of the transformed equations using integrating
factors, setting:
\begin{equation}
U=\e^{-\Omega^2 t}\tilde{U}, \quad
V=\e^{-\epsilon\Omega^2 t}\tilde{V},
\end{equation}
such that now
\begin{equation}
\partial_t\tilde U=\e^{\Omega^2 t}\cF[f(u,v)], \quad
\partial_t\tilde V=\e^{\epsilon\Omega^2 t}\cF[g(u,v)].
\label{eq:ode}
\end{equation}

At this point, in practical terms, we discretize the spatial
domain, considering $N_x$ and $N_y$ equispaced points in the $x$
and $y$ directions. Then we utilize the discrete FFT so equation
(\ref{eq:ode}) becomes a system of ODEs parameterized by the
Fourier modes (distinguished by a couple of indices $ij$) so
\begin{equation}
\partial_t\tilde U_{ij}=\e^{\Omega^2_{ij} t}\cF[f(u_{ij},v_{ij})], \quad
\partial_t\tilde V_{ij}=\e^{\epsilon\Omega^2_{ij} t}\cF[g(u_{ij},v_{ij})],
\label{eq:ode_discretized}
\end{equation}
where $u_{ij}=u(x_i,y_j)$, $v_{ij}=v(x_i,y_j)$ and
$\Omega^2_{ij}=\omega_{x}^2(i)+\omega_{y}^2(j)$. Periodic boundary
conditions are implicitly set at the extremes of the spatial
domain. Henceforth we will suppress the $ij$ indices and consider
this discretization understood. In practice one takes, say,
$N_x=N_y=N=128$ and utilize $128\times 128$ Fourier modes.

The spatial derivatives have now disappeared, along with the
stiffness that they had introduced, and the resulting ODEs are now
simple to solve with an explicit solver, say, Runge-Kutta. There
are issues that arise at this stage, for instance, the fixed
points of these equations are different from those of the original
ones, aliasing could be an issue and, if we use Runge-Kutta or
some other scheme, how does the local truncation error depend upon
$\Omega_{ij}$?
 We leave these issues aside until
later in the article.

Adopting the standard notation for order $M$
Runge-Kutta methods, with time step $\Delta t$, that is to advance
from $t_n=n\Delta t$ to $t_{n+1}=(n+1)\Delta t$ for an ODE
$y_t=f(t,y)$: $y_{n+1}=y_n+\sum_{i=1}^M c_i k_i$ and each $k_i$ is
$k_i=\Delta t f(t_n+a_i\Delta t, y_n+\sum_{j=1}^{i-1}b_{ij}k_j)$.
The $a_i, c_i, b_{ij}$ are given by the appropriate Butcher array;
an almost infinite number of different schemes exist. Two popular
ones are given in Numerical Recipes \cite{press92a}: classical
fourth order and the Cash-Karp embedded scheme \cite{cash90a}, we
utilize these in our numerics.

For our purposes we apply the general explicit Runge-Kutta formula
to the equations (\ref{eq:ode_discretized}) for $\tilde U_{ij}$
and $\tilde V_{ij}$. 
Notationally, we denote $\mu_i$ and $\nu_i$ to be the $k$'s
associated with the $\tilde U$ and $\tilde V$ equations
respectively. The right-hand sides have slightly unsettling
exponential terms in $t$ and it is convenient to set replacement
variables as
\begin{equation}
\tilde \mu_i = \mu_i \e^{-\Omega^2 t_n} , \quad
\tilde \nu_i = \nu_i \e^{-\epsilon\Omega^2 t_n}.
\end{equation}
We write the formulae out for $U$ and $V$ since it is simpler to
just work with the transforms of the physical variables rather
than the physical variables themselves. The upshot is that
with an $M$-stage Runge-Kutta scheme,
\begin{equation}
U_{n+1}=  e^{-\Omega^2\Delta t}\left[U_n +\sum_{i=1}^M c_i \tilde \mu_i\right], \quad
V_{n+1}=  e^{-\epsilon\Omega^2\Delta t}\left[V_n +\sum_{i=1}^M c_i \tilde \nu_i\right],
\end{equation}
where the modified $\tilde \mu_i$ and $\tilde \nu_i$ terms are
\begin{eqnarray}
\tilde \mu_i & = & \e^{\Omega^2 a_i\Delta t}\Delta t
\cF\left\{f\left[\cF^{-1}\left(U_{n+a_i}\right),
\cF^{-1}\left(V_{n+a_i}\right)\right]\right\}, \nonumber \\
\tilde \nu_i & =& \e^{\epsilon\Omega^2 a_i\Delta t}\Delta t
\cF\left\{g\left[\cF^{-1}\left(U_{n+a_i}\right),
\cF^{-1}\left(V_{n+a_i}\right)\right]\right\},
\end{eqnarray}
and the values of $U$ and $V$ at the intermediate steps are
\begin{equation}
U_{n+a_i} = e^{-\Omega^2 a_i\Delta t}\left[U_n +
\sum_{j=1}^{i-1} b_{ij}\tilde \mu_j \right], \quad
V_{n+a_i} = e^{-\epsilon\Omega^2 a_i\Delta t}\left[V_n  +
\sum_{j=1}^{i-1} b_{ij}\tilde \nu_j\right].
\end{equation}
That is, one works entirely in the spectral domain and one inverts
a transform to recover $u$ and $v$. Clearly a fair amount of
Fourier transforming to and from, is involved and this is the
primary numerical cost. Fortunately Matlab has simple to use
multi-dimensional Fast Fourier Transform (FFT) routines, ({\tt
fft, ifft, fft2, ifft2}) and many routines are available in
Fortran (or C).

The essential point is that by removing the stiffness one can use
explicit high-order timesolvers and rapidly and accurately move
forwards in time, this is vastly superior to using implicit
schemes particularly in higher dimensions (see
Appendix~\ref{app:adi}). There are some slight deficiencies in the
method that can be removed using more recent ideas and we shall
return to the theory in section \ref{sec:etd}.

\section{Illustrative examples}
We choose a range of illustrative examples that are of current and
recurring interest, and which cover pitfalls and natural questions
that arise.

\subsection{One dimensional models}

\subsubsection{Fisher's equation: Speed selection}
Fisher's equation
\begin{equation}
    u_t=u_{xx}+u(1-u),\quad \vert x\vert<L
\label{eqn:fisher1D}
\end{equation}
provides a nice demonstration; there is a detail that is worth
investigating: the speed selection associated with exponential
decay of the initial condition. Numerically this has been an issue
for other approaches such as moving mesh schemes (\cite{qiu98a})
with some authors recommending that such schemes be used with
caution upon this type of problem (\cite{li98a}).

\begin{figure}
\leavevmode
\centering
\includegraphics[height=8.5cm]{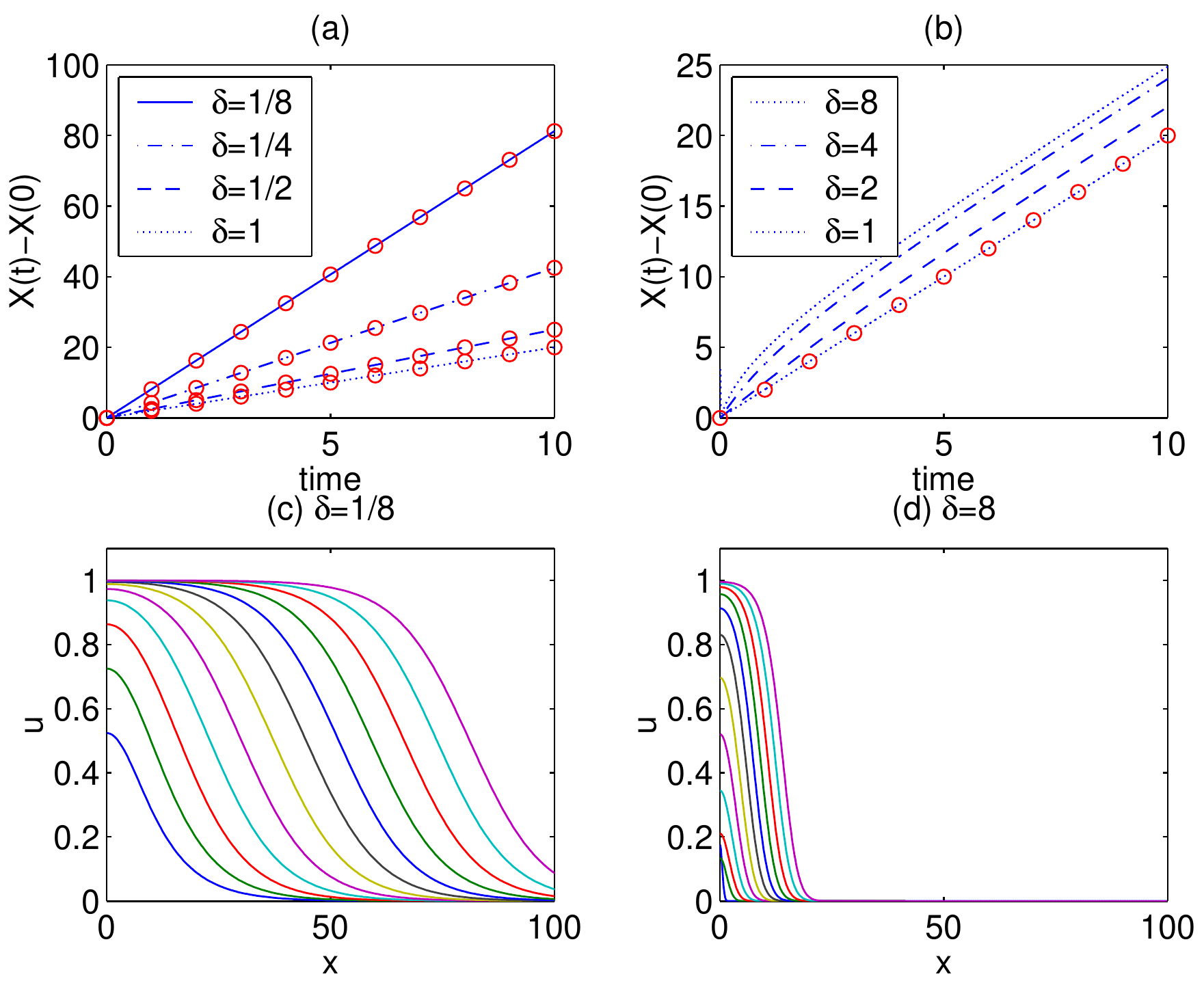}
\vspace*{2mm}
\caption{Panels (a) and (b) show front positions versus time for
various values of $\delta$. In panel (a) these are shown for
$\delta\le 1$. Panel (b) shows front positions for $\delta\ge 1$.
Circles show $ct$, $c$ taken from (\ref{eq:c}). Panels (c) and (d)
show profiles of $u$ at unit time intervals until $t=10$. }
\label{fig:fig1}
\end{figure}

We utilize an initial condition
\[
   u(x,0)=\frac{1}{2\cosh\delta x}
\]
that has exponential decay $\exp(-\delta |x|)$ as $\vert
x\vert\rightarrow\infty$. Theoretically, one expects travelling
waves to develop from such an initial condition on an infinite
domain, which we truncate at some large, but finite, value, say,
$L\sim 150$, and what is particularly interesting is that the
system then selects the constant velocity at which the developed
fronts propagate, $c$, and the velocity is a function of the decay
rate of the initial condition:
\begin{equation}
    c=2\quad {\rm for}\ \delta>1, \quad c=\delta+\frac{1}{\delta}\quad {\rm otherwise}.
\label{eq:c}
\end{equation}
We easily extract the velocity from the simulations, figure
\ref{fig:fig1}, and evidently we recover these theoretical values.
It is essential that $L$ is taken to be large enough that the
initial condition is effectively zero at $L$, here $L\sim 150$.
The figure shows the front position, $X(t)$, here taken as the
point where $u(X,t)=10^{-4}$, versus time, the slope gives the
velocity; for convenience $ct$ is also shown, note the offset in
panel (b) is not relevant as it is the slope that concerns us.

The simulations in Fortran 77  are performed virtually
instantaneously, with the Matlab simulations taking a few seconds.

\subsubsection{Gray-Scott: Pulse splitting}
Another very interesting feature of many reaction diffusion
equations is pulse splitting or shedding; a propagating pulse is
unstable, and the unstable eigensolutions lag behind the pulse
causing a daughter pulse to break off. This is particularly
pronounced in the Gray-Scott equations, fortunately they have a
pleasant and simple nonlinearity in the reaction terms that makes
them amenable to analytical approaches
\cite{doelman97a,reynolds97a}. They also form a tough test
upon any numerical scheme as the splitting events and subsequent
structure must be captured correctly both in space and time.

The equations ($v$ is the activator and $u$ the inhibitor) are:
\begin{equation}
u_t=u_{xx}-uv^2+A(1-u),\quad
v_t=\epsilon v_{xx}+u v^2 -Bv .
\label{eqn:gray1D}
\end{equation}

\begin{figure}
\leavevmode
\centering
\includegraphics[height=5cm]{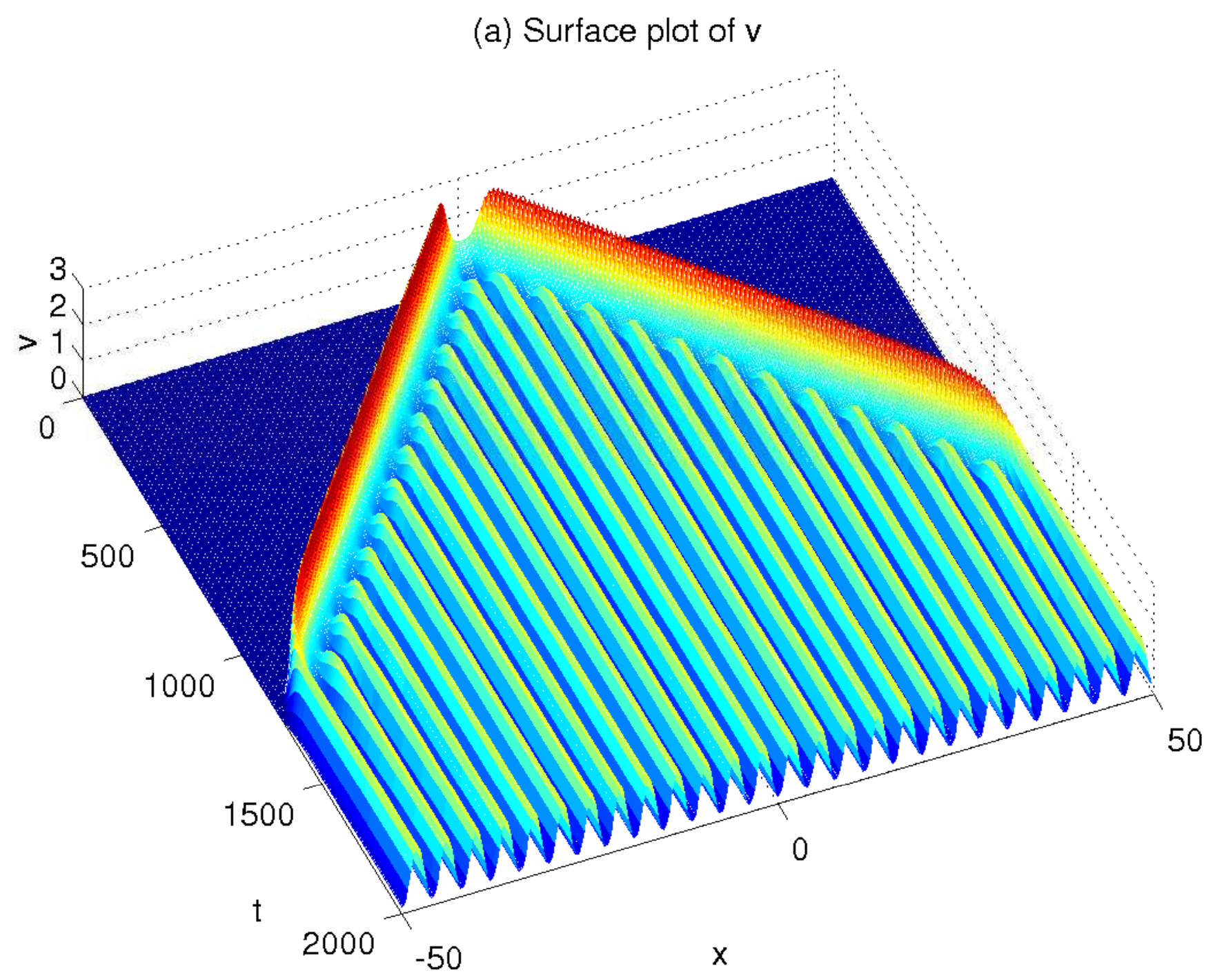}
\includegraphics[height=4.5cm]{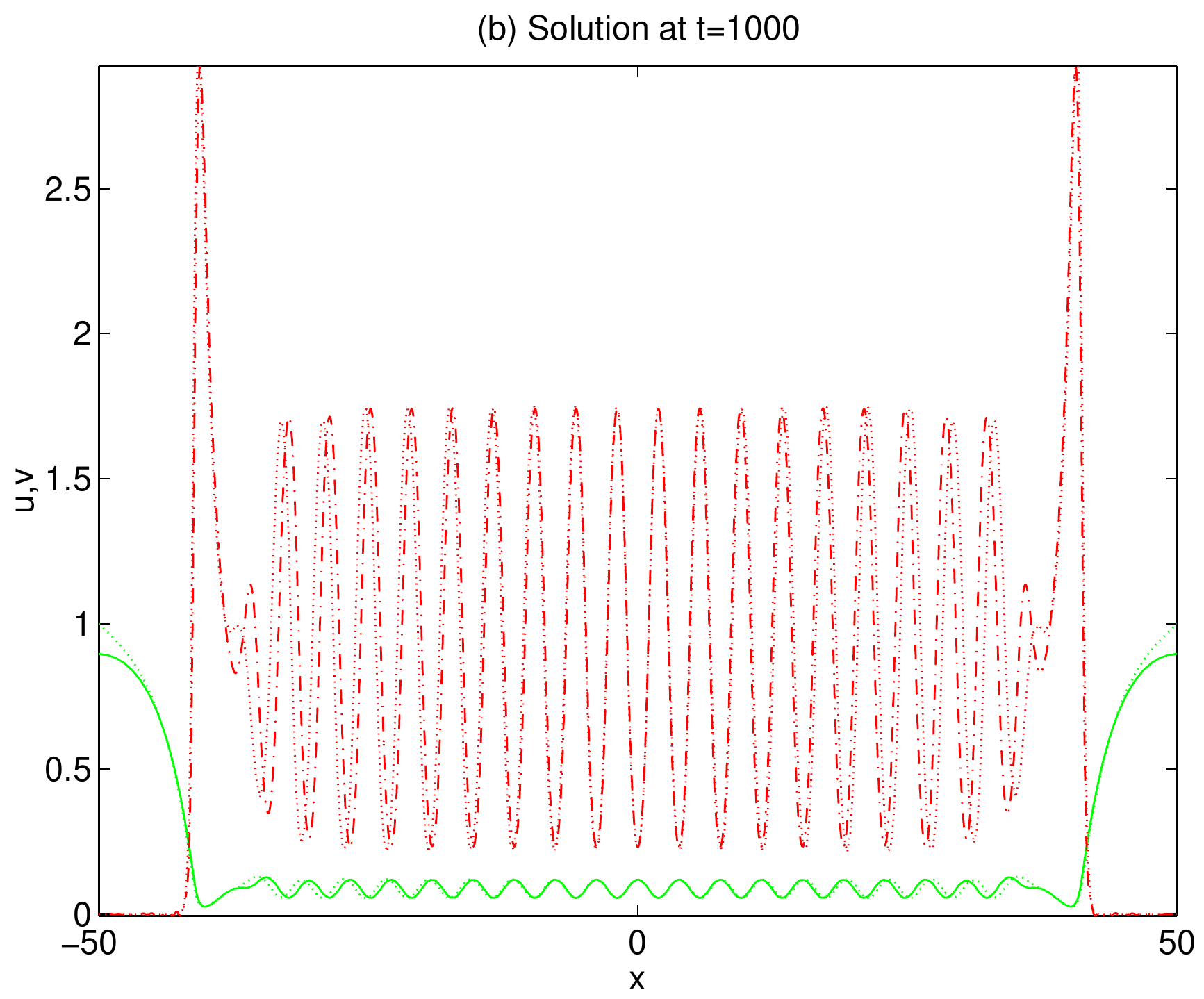}
\vspace*{2mm}
\caption{Panel (a) shows the outward propagating pulses and the
shedding phenomena for $v$ out to $t=2000$. Panel (b) shows both
$u$ (solid) and $v$ (dot-dashed) at $t=1000$, the dotted lines
also shown come from the adaptive scheme [Blom \& Zegeling 1994].
The parameters chosen here are $a=9, b=0.4, \epsilon=0.01$ where
$A=\epsilon a, B=\epsilon^{1/3}b$.} \label{fig:gscott1D}
\end{figure}

A couple of illustrative plots are given in figure
\ref{fig:gscott1D}, these have initial conditions
\[
u=1-\frac{1}{2}\sin^{100}(\pi(x-L)/2L), \quad
v= \frac{1}{4}\sin^{100}(\pi(y-L)/2L)
\]
where we choose the half domain length, $L$, to be $50$. These are
chosen to replicate a figure from \cite{doelman97a}; notably the
simulations differ as the boundary conditions here at $\pm L$ are
periodic, and eventually a steady spatially periodic state
emerges. The simulations in \cite{doelman97a}, and reproduced in
panel (b) of figure \ref{fig:gscott1D}, utilize Dirichlet, that
is, fixed values of $u,v$, conditions and hence there is a minor
discrepancy close to the edges of the domain, this is most
noticeable in $u$. The Matlab file for this computation is given
in Appendix~\ref{app:code}, and it takes 12 seconds to run on a
1GHz Pentium 3 Dell Laptop running Linux. Comparative computation
times are probably meaningless as computational power will ever
increase, our only point being that these computations are
relatively fast versus competitors even with modest computing
facilities available to all/ most undergraduates. The adaptive
scheme takes a few minutes depending upon the number of grid
points utilized, typically 500 points.

\subsubsection{Autocatalysis: Oscillatory fronts}
\begin{figure}
\leavevmode
\centering
\includegraphics[height=5cm]{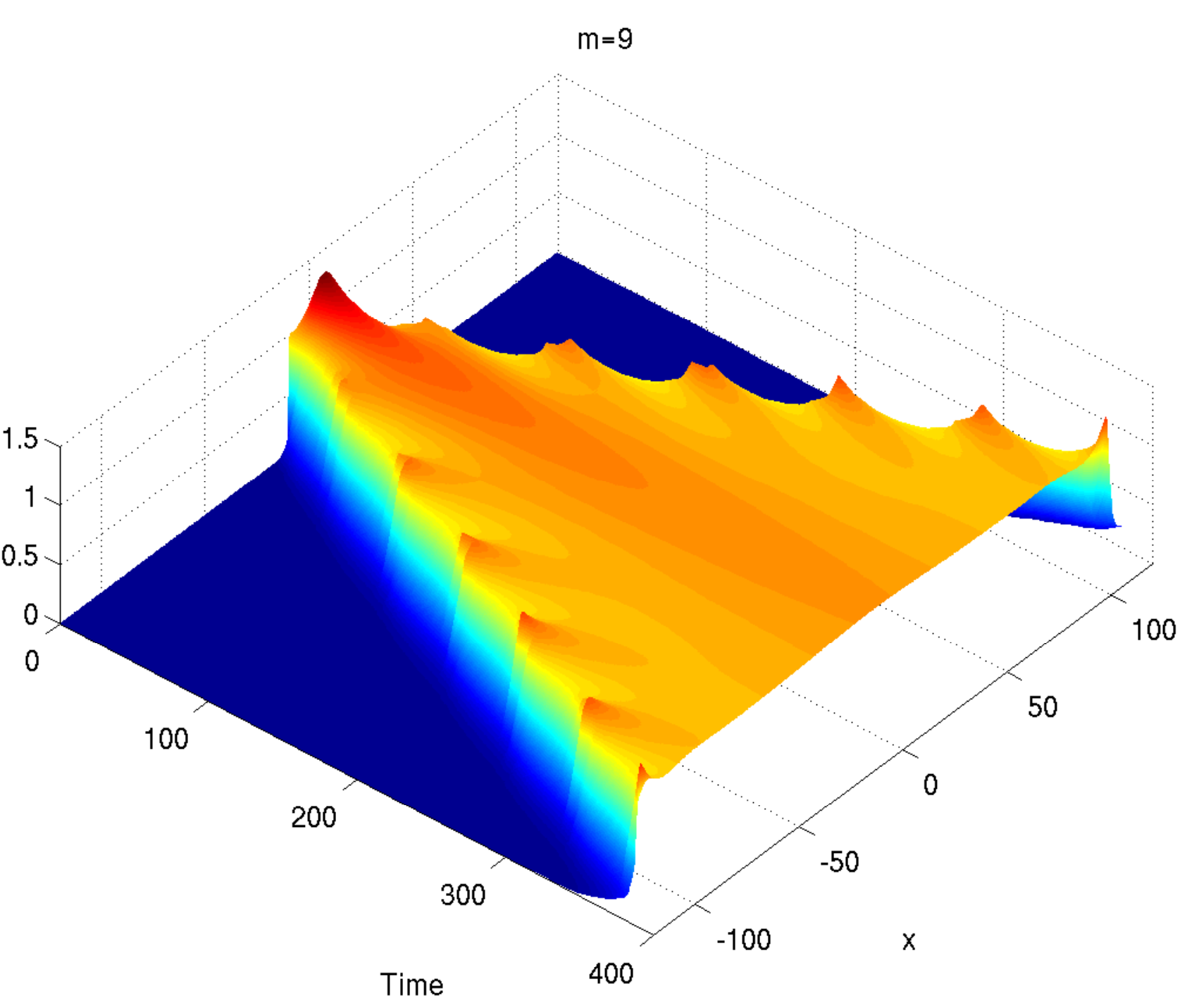}
\includegraphics[height=5cm]{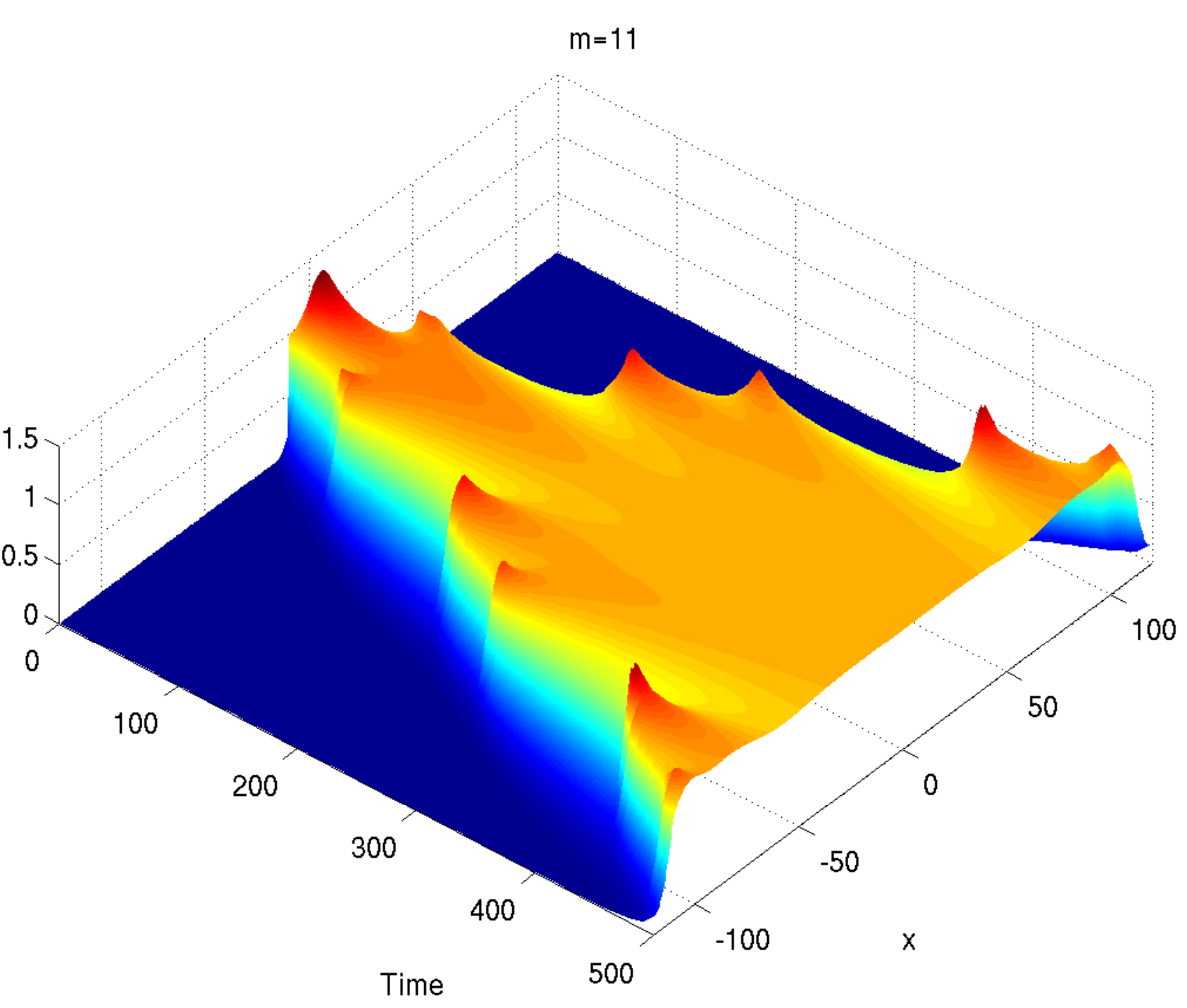}
\vspace*{2mm}
\caption{Typical results for the autocatalytic model with
$\epsilon=0.1$ and $m=9$ and $m=11$. } \label{fig:fig2}
\end{figure}

Many reaction-diffusion equations arise in combustion theory, or
in related chemical models. One such model, in non-dimensional
terms, is
\begin{equation}
u_t = u_{xx} + vf(u),\quad v_t=\epsilon v_{xx}-vf(u),
\label{eqn:auto}
\end{equation}
where
\begin{equation}
f(u)=\cases{u^{m}, \ \ u\ge 0,\cr 0,\ \ u<0}
\label{eqn:fu}
\end{equation}
and $\epsilon$ is the inverse of the Lewis
number (the ratio of diffusion rates). It arises when two chemical
species $U$ and $V$ react such that $mU+V\rightarrow (m+1)U$; the
two species have different diffusivities, their ratio being $\epsilon$
($\epsilon<1$ is the regime of interest). What is particularly
interesting in this model is that steady travelling waves occur
for low values of $m$, their speed is a function of $m, \epsilon$ and
the fronts steepen dramatically for large $m$. In fact, in full
nonlinear simulations as $m$ increases a Hopf bifurcation occurs,
and as it increases yet further one gets chaotic behaviour at the
wave front. This behaviour is detailed in \cite{balmforth99a},
\cite{Metcalf:1994} and similar features arise in combustion
models, see for instance \cite{BaylissMat:1989} where $f(u)$ is
replaced by exponential Arrhenius reaction terms. For our
verification purposes we compared with the computations of \cite{balmforth99a},
 and initiated the computations with
\[
u=\frac{1}{2}\left(1+\tanh(10(10-|y|))\right) , \quad
v=1-\frac{1}{4}\left(1+\tanh(10(10-|y|))\right)
\]
that is, a sharp localized disturbance and obtained perfect
agreement even for extremely steep fronts. Figure \ref{fig:fig2}
shows typical results with periodic fluctuations at $\epsilon=0.1$
and $m=9$ with more extreme behaviour with $m=11$. It is notable
that all the delicate behaviour, rocking fronts and transitions to
apparently chaotic behaviour, is accurately captured, together
with the very steep fronts; this is a challenge for any numerical
scheme, a minor aside is that our Fortran 77 code had difficulties
with this computation until we imposed symmetry conditions at the
end of every timestep, or zeroed the imaginary part of the
inverted result.

This leads us to an algorithmic detail: in the Matlab codes we use
separate variables $u$ and $v$ and transform, and invert, each
independently and just use the real parts of each inverse -
automatically zeroing the imaginary parts and thereby preventing
rounding errors mounting up over time. However, this is actually a
bit wasteful as one could combine $u$ and $v$ to be the real and
imaginary parts of a single transform variable and just halve the
amount of work; the Fortran codes use the latter approach. This
only seems to cause problems in this particular example as the
extreme powers $11$ magnify rounding errors and, for the Fortran
code, the minor modification described above is required to
maintain accuracy.

\subsection{Two dimensional examples}
It is in higher dimensions that the ideas presented here really
become of serious value. We choose to illustrate the numerical
algorithms using a couple of non-trivial examples from the
reaction-diffusion equation literature. Stable, labyrinthine
patterns arising in FitzHugh-Nagumo type reaction diffusion
equations \cite{hagberg94a}, and pulse splitting from the
Gray-Scott equations \cite{muratov01a}.

\subsubsection{Gray-Scott}

If we have radial symmetry then many of the 1D schemes can be
utilized in radial coordinates, an adaptive code \cite{blom94a} is
particularly convenient as it takes advantage of the
\cite{skeel90a} discretization that automatically incorporates the
coordinate singularity. We can therefore check the numerical
simulations of the 2D spectral code versus these. The Gray-Scott
model has delicate features that are not easy for a numerical
scheme to extract, and it is susceptible to small perturbations
generating instabilities.

\begin{figure}
\leavevmode
\centering
\includegraphics[height=8cm]{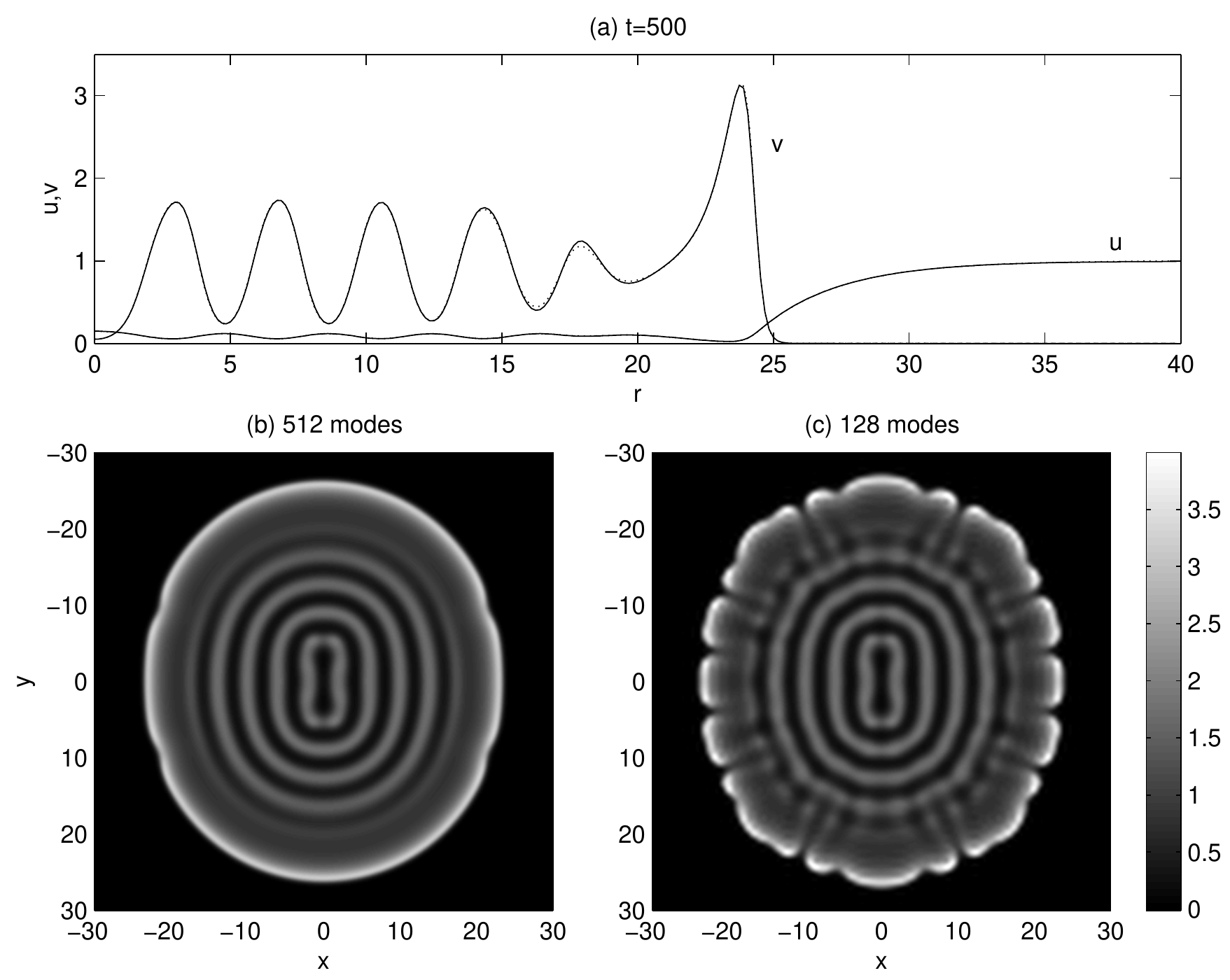}
\vspace*{2mm}
\caption{Panel (a) shows the outward propagating rim and the
shedding phenomena for $u,v$, at $t=500$ for an axi-symmetric
initial condition, along a radial line; the barely visible dotted
lines also shown come from the adaptive scheme [Blom \& Zegeling
1994]. Panels (b) and (c) shows a planview of $v$ at $t=500$
starting from non-axisymmetric initial data, the two panels show
$512$ and $128$ modes respectively (same parameters as figure
\ref{fig:gscott1D}). For a better display, resolution in panel (c)
was improved using Fourier cardinal function interpolation (see
Appendix~\ref{app:cardinal}). } \label{fig:gs2D}
\end{figure}

As noted earlier, the Gray-Scott equations have provided an
interesting test bed for theoreticians exploring pulse splitting
and so-called auto-solitons and their stability
\cite{doelman97a,reynolds97a,muratov01a,muratov02a}; it is also
notable that the model also arises in biological contexts
\cite{davidson97a}. Different authors prefer different rescalings
according to the physics/biology that they wish to emphasise, we
shall not enter that debate here. The equations we use are
\begin{equation}
u_t=u_{xx}+u_{yy}-uv^2+A(1-u),\quad v_t=
\epsilon [v_{xx}+v_{yy}]+u v^2 -Bv .
\label{eqn:gray2D}
\end{equation}

The two-dimensional analogue of the pulse-splitting events of
figure \ref{fig:gscott1D} are shown in figure \ref{fig:gs2D};
using axisymmetric initial conditions for panel (a):
\[
    u=1-\frac{1}{2}\exp( -r^2/20), \quad v=\frac{1}{4}\exp(- r^2/20)
\]
where $r^2=x^2+y^2$, allows us to verify the two-dimensional
computations in a non-trivial way as the system is highly
unstable. Although not shown, it is particularly striking how
perfect axisymmetry is retained by these axisymmetric
computations. Altering the initial conditions to break the
axisymmetry to the above, but with $r^2= x^2/2+y^2$ leads to the
oval pattern of alternating high and low concentrations shown in
figure \ref{fig:gscott1D} panels (b) and (c). What is particularly
notable is that using fewer modes leads to an attractive, but
evidently erroneous, pattern.

\subsubsection{Labyrinthine Patterns}
\begin{figure}
\centering
\includegraphics[height=8.5cm]{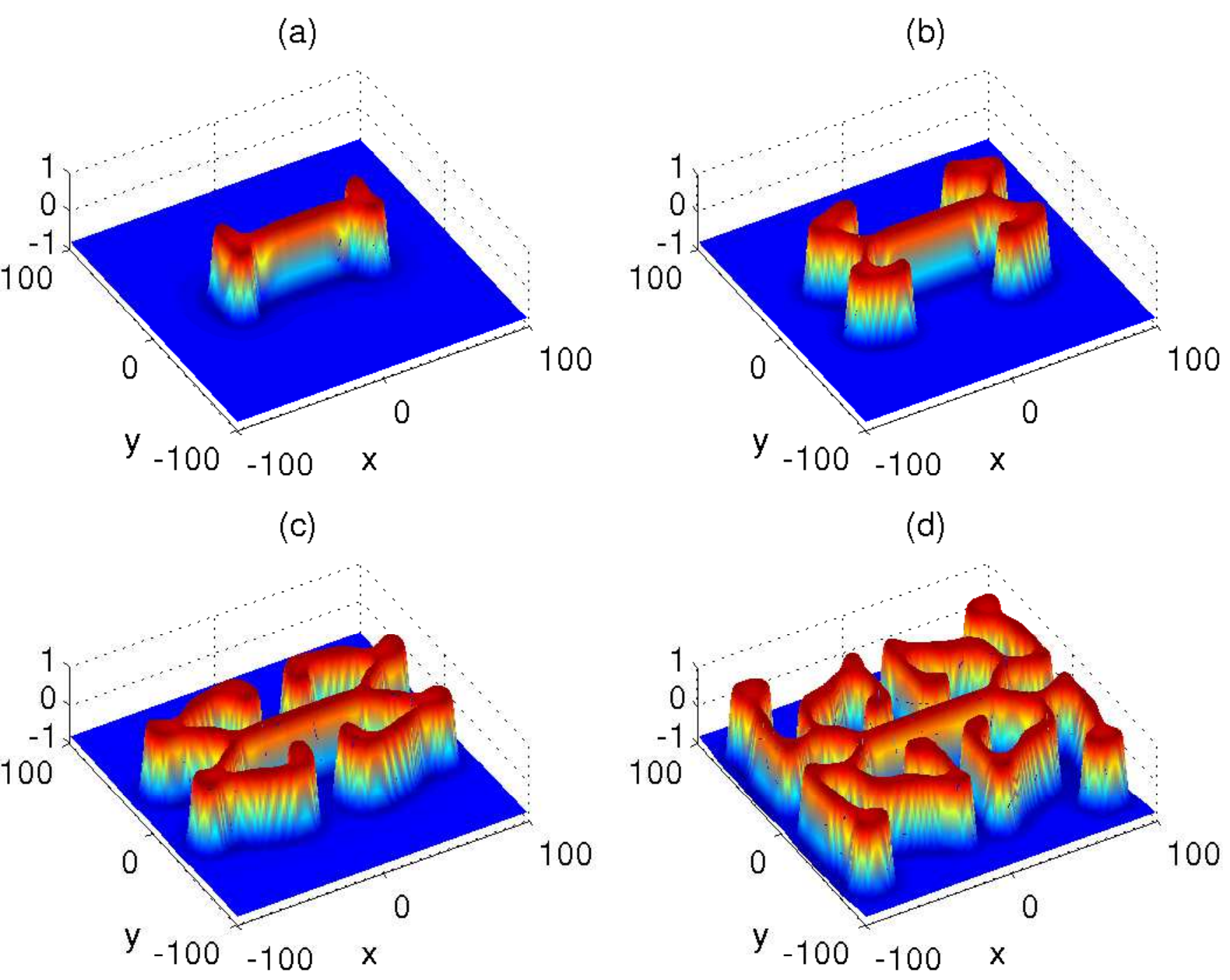}
\includegraphics[height=4cm]{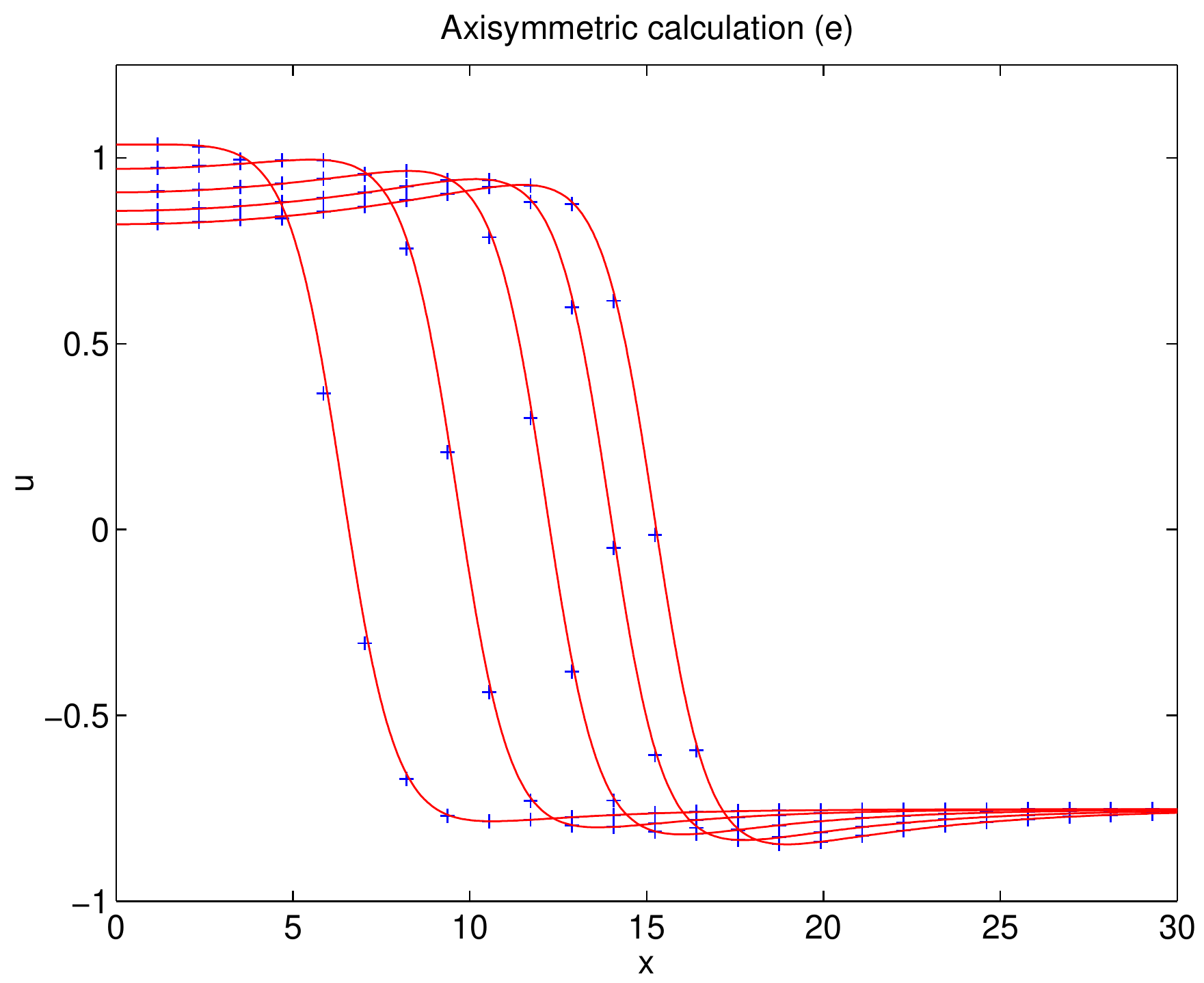}
\includegraphics[height=4cm]{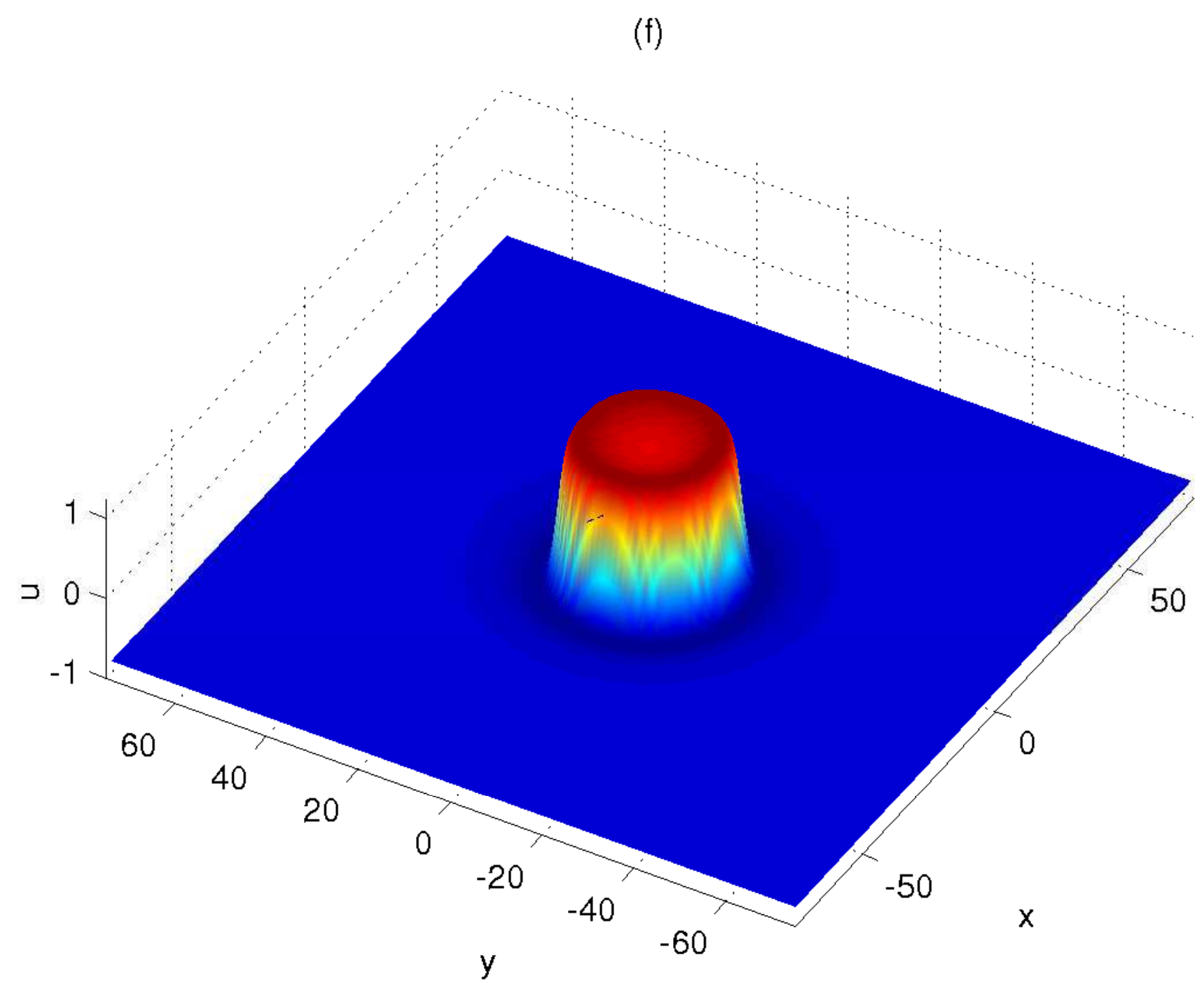}
\vspace*{2.5mm}
\caption{Panels (a) to (d) show the emerging labyrinthine pattern
for $u$ at times $t=200,400,600,1000$. Panels (e) and (f) an
axisymmetric computation, same initial conditions as equation
(\ref{eq:initialfitz}) bar that $0.01y^2\rightarrow y^2$, at
$t=10,20,30,40,50$, in (e) the solid lines come from the adaptive
scheme of [Blom \& Zegeling 1994] and crosses from the spectral
method. The parameters chosen here are $a_0=-0.1, a_1=2,
\epsilon=0.05, \delta=4$. } \label{fig:labyrinthine}
\end{figure}

A striking and interesting group of patterns that emerge in models
of catalytic reactions are growing labyrinthine patterns
\cite{hagberg94a,meron01a}. Starting from a non-axisymmetric
initial condition strongly curved portions move more rapidly and
the pattern lengthens. Regions of high concentrations repel and,
hence from the periodicity of the domain, the patterns turn inward
until an equilibrium is reached. An illustrative simulation is
given in figure \ref{fig:labyrinthine}. The computation utilizes
$128\times 128$ Fourier modes on a $200\times 200$ grid. Doubling
the number of Fourier modes makes no discernable difference; more
detailed numerical error discussions are in a later section.

The governing equations are that
\begin{equation}
u_t=u-u^3-v+\nabla^2 u, \quad
v_t=\delta(u-a_1 v-a_0)+\epsilon\nabla^2 v
\label{eqn:labyrinthe2D}
\end{equation}
where $u, v$ represent activator and inhibitors. The parameters
$a_0, a_1$ , $\epsilon, \delta$ lead one from one regime to
another, see \cite{hagberg94a} for details. We begin the
simulation in figure \ref{fig:labyrinthine} from initial
conditions
\begin{equation}
u=a_1 v_-+a_0- 4a_1 v_-e^{-0.1(x^2+0.01 y^2)}, \quad
v=v_--2v_-e^{-0.1(x^2+0.01 y^2)}.
\label{eq:initialfitz}
\end{equation}
This being an elliptical mound of chemical concentrations of
sufficient magnitude to trigger the reaction. Here
$v_-=(u_--a_0)/a_1$ is found from, $u_-$, which is the smallest
real root of the cubic $a_1 u^3+u(1-a_1)-a_0=0$. This $(u_-,v_-)$
state being stable. As one notes from the figure, the evolution is
non-trivial and the edges of the concentrations are sharp and
steep, all features that test the robustness of the scheme.
Verification follows from comparison with the adaptive scheme for
axisymmetry, and since a stationary state emerges one can also
generate the final shape from a boundary value problem; all
computations agree. Figure \ref{fig:labyrinthine} only shows the
evolution of $u$, that of $v$ is qualitatively similar.

It is worth noting that other behaviours are possible for these
equations in other parameter regimes than those chosen here.

\subsection{Refinements}
All of the figures shown are generated using a fixed time step,
and the classical standard fourth order Runge-Kutta scheme; this
was $0.1$ in all cases except the autocatalytic problem with large
$m$ where we took a timestep of $0.02$. This is deliberate to
demonstrates that sophisticated algorithms are not vital, however
it is not pleasant to have no error control or indeed no idea of
how accurate the solution actually is at each time step. An overly
enthusiastically large choice for the time step could lead to
numerical instabilities and accumulated error. The results we
present are all generated using Matlab; the Matlab codes that we
present in the appendices are efficient as teaching tools, and, to
a certain extent as research tools; the high level language gives
short and easily understandable code. However, traditional
languages such as Fortran or C, C++ will often run much faster, at
least versus interpreted Matlab code, and to complement the Matlab
codes we also provide the source codes in Fortran 77.

\subsubsection{Adaptive time stepping}

Adaptive time stepping based upon embedded 5th order Runge-Kutta
schemes in the usual manner, see for instance
\cite{cash90a,press92a}, is easily implemented and we do so. The
user can replace these weights with their favourite scheme
\cite{dormand80a}, say, but this will change little in practice.
These codes incorporate an error tolerance and utilize local
extrapolation, these codes often settle to using surprisingly
large time steps $dt\sim 0.5$ or larger for even moderate error
tolerances (relative error $\sim 10^{-4}$) and this  further
speeds the computations.

\subsubsection{Exponential time differencing}
\label{sec:etd}

Integrating factor ideas are not new, there are actually several
non-appealing aspects of the approach. On the basis of ``truth in
advertising'' we must reveal them. A philosophically unpleasant
feature is that one notes that the fixed points of equations
(\ref{eq:ode}) are not the same as those of the original
untampered equations (\ref{eq:original_odeU}),
(\ref{eq:original_odeV}). We are unaware of any circumstance in
which this has led to any problems, but it is not nice. A more
serious fact that counts in its disfavour is that the local
truncation error, when $\Omega^2 \Delta t\ll 1$, for the time
stepping schemes we use are, for say 4th order Runge-Kutta,
$O(\Omega^2 \Delta t)^5$. This is apparently disastrous as
$\Omega$ can be large and we now have fourth order accuracy in
time, but with a large pre-multiplicative factor. However, it is
important to recall that the local truncation error involves
expanding $\exp(-\Omega^2 \Delta t)$ terms which are, practically,
exponentially small for large $\Omega$ relative to $\Delta t$.
Nonetheless we are clearly picking up an additional contribution
to the numerical error for moderate values of $\Omega$. This is
surmountable, basically one designs a time stepping scheme that
correctly incorporates the exponential behaviour, a recent
article, \cite{cox02a}, derives several exponential Runge-Kutta
schemes; this is an active area of current research with
modifications of their scheme by \cite{kassam03a} to overcome a
numerical instability and by \cite{krogstad05a} generating a
scheme with smaller local truncation error and better stability
properties. Krogstad also notes the very interesting link with the
commutator-free Lie group methods of \cite{munthe-kaas99a} and
undoubtedly this area will develop further.

In essence the exponential time differencing idea applied here for
$U$, involves utilizing the integrating factor $\exp(\Omega^2 t)$
and multiplying equation (\ref{eq:original_odeU}) through by it
and then we integrate over a time step to obtain:
\[
  U_{n+1}=U_n e^{L \Delta t} +e^{L \Delta
  t}\int_{0}^{\Delta t} e^{-L\tau} N_u(u(t_n+\tau),v(t_n+\tau),t_n+\tau) d\tau
\]
where we have rewritten equation (\ref{eq:original_odeU}) as
\[
   U_t=Lu +N_u(u,v)
\]
that is, with a Linear piece (here $-\Omega^2$) and a Nonlinear
piece (here the Fourier transform of the nonlinear reaction
terms); this is the notation used in the relevant literature. The
interesting departure, and distinguishing feature, from the
standard integrating factor method is that one then approximates
the integral and the truncation error is then independent of
$\Omega^2$. The article by \cite{cox02a} contains various
approximations to the integral and numerical comparisons of
methods.

It is important to bring these ``state-of-the-art'' solvers into
the more applied domain and enable other researchers to take
advantage of them. We utilize the fourth order Runge-Kutta-like
scheme of Krogstad,
\begin{eqnarray*}
U_{n+1}& = & e^{L\Delta t}U_n+\Delta t[4\phi_2(L\Delta t)-3\phi_1(L\Delta t)+
\phi_0(L\Delta t)] N_u(U_n,V_n,t_n)+ \\
& & \hspace{0.5cm} 2\Delta t[\phi_1(L\Delta t)-2\phi_2(L\Delta t)]
N_u(\mu_2,\nu_2,t_n+\Delta t/2)+\\
& & \hspace{0.5cm} 2\Delta t[\phi_1(L\Delta t)-2\phi_2(L\Delta t)] N_u(\mu_3,\nu_3,t_n+\Delta t/2)+ \\
& & \hspace{0.5cm} \Delta t[ 4\phi_2(L\Delta t)-\phi_1(L\Delta t)]N_u(\mu_4,\nu_4,t_n+\Delta t)
\end{eqnarray*}
with the stages $\mu_i$ as
\begin{eqnarray*}
\mu_2 & = & e^{L\Delta t/2}U_n+(\Delta t/2)\phi_0(L\Delta t/2) N_u(U_n,V_n,t_n)\\
& & \\
\mu_3 & = & e^{L\Delta t/2}U_n+(\Delta t/2)\left[ \phi_0(L\Delta t/2)-
2\phi_1(L\Delta t/2)\right] N_u(U_n,V_n,t_n)+\\
& & \hspace{0.5cm} \Delta t \phi_1(L\Delta t/2)N_u(\mu_2,\nu_2,t_n+\Delta t/2) \\
& & \\
\mu_4 & = & e^{L\Delta t}U_n+\Delta t\left[ \phi_0(L\Delta t)-
2\phi_1(L\Delta t)\right] N_u(U_n,V_n,t_n)+ \\
& & \hspace{0.5cm} 2\Delta t \phi_1(L\Delta t)N_u(\mu_3,\nu_3,t_n+\Delta t).
\end{eqnarray*}

The functions $\phi_i$ are defined as
\[
    \phi_0(z)=\frac{e^z-1}{z}, \quad
    \phi_1(z)=\frac{e^z-1-z}{z^2} \quad
    \phi_2(z)=\frac{e^z-1-z-z^2/2}{z^3}
\]
and these are precisely the terms that emerge naturally in the Lie
group methods \cite{munthe-kaas99a}. The original Cox \& Matthews
scheme involves a split-step and has marginally worse error and
stability properties. There are also slight problems associated
with capturing the behaviour of $\phi_i(z)$ uniformly as
$z\rightarrow 0$ and \cite{kassam03a} suggest a remedy; one uses
an integral in the complex plane \cite{higham96a}, and we also use
this approach in our algorithms. For brevity we have presented the
scheme for $U$ alone, the $V$ equations follow in a similar
fashion.

\begin{figure}
\leavevmode
\centering
\includegraphics[height=4.5cm]{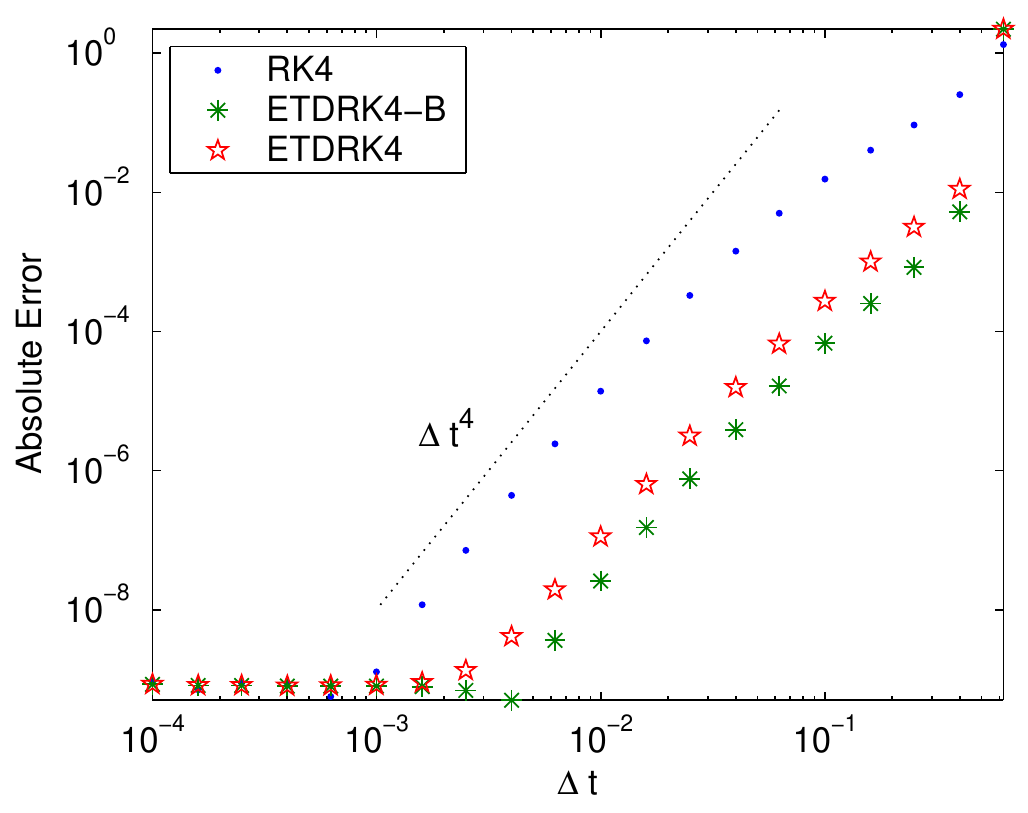}
\vspace*{2mm}
\caption{The $1$-D Gray-Scott equation~(\ref{eqn:gray1D}) is
solved using different time steps $\Delta t$ with the same
parameters values of figure (\ref{fig:gscott1D}). The final
solution obtained at $t=200$ is compared, for each method (RK4,
ETDRK4 and ETDRK4-B), with a gold-standard run (computed with
ETDRK4-B and $\Delta t=10^{-5}$); the maximum absolute errors are
displayed as a function of the time step.}
\label{fig:error_comparison_gray}
\end{figure}

We label this as a fourth order Exponential Time Differencing
Runge-Kutta (ETDRK4-B) scheme to distinguish it from a standard
Runge-Kutta scheme and to be consistent with the notation of
\cite{cox02a,kassam03a,krogstad05a}. It is also worth noting that
various other exponential {R}unge-Kutta schemes have been
developed by other authors, see \cite{berghe00a} and the
references therein to overcome this difficulty arising in other
contexts. \nocite{cox02a,berghe00a}

A reasonably large numerical overhead is involved in setting up
the fourth-order scheme using the device suggested by
\cite{kassam03a}; if one utilizes an adaptive scheme in time then
this overhead must be regularly recomputed and this then becomes
expensive, hence we do not adapt the ETDRK4 schemes in time.

A numerical $1$-D comparison of the ETDRK4 scheme of
(\cite{cox02a}), the improved ETDRK4-B (\cite{krogstad05a}), and
the more standard RK4 scheme is shown in figure
(\ref{fig:error_comparison_gray}). It is noticeable that ETDRK4-B
proved to give the smaller error; both it and the ETDRK4 scheme
provide an order of magnitude improvement over the RK4 scheme for
larger timesteps, rewarding the extra programming effort. The
errors all scale with the expected $\Delta t^4$ scaling.

Although from a practical point of view all of these schemes are
explicit and so are much better (faster, accurate) than the
implicit, or semi-implicit schemes often used for these equations.
This must be the main message to be taken from this article.

\begin{figure}
\leavevmode
\centering
\includegraphics[height=4.5cm]{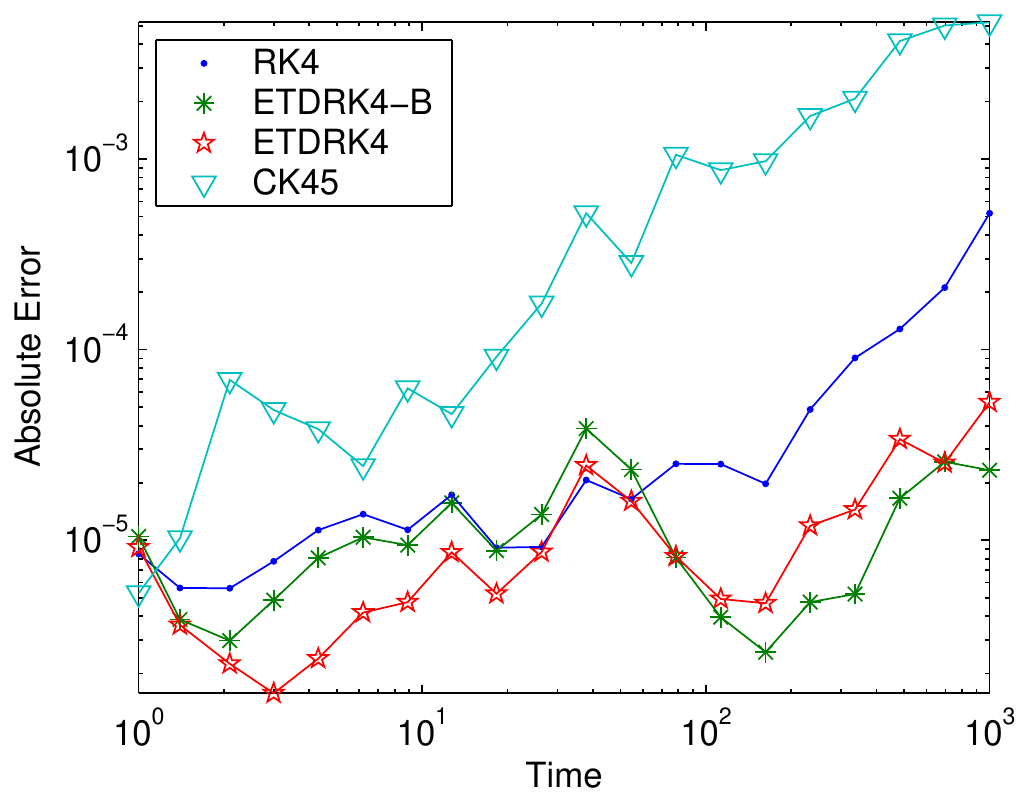}
\vspace*{2mm}
\caption{
Absolute errors of the computational methods with respect to a
gold-standard run obtained with ETDRK4-B using $\Delta t=0.01$.
The $2$-D equation being solved is (\ref{eqn:labyrinthe2D}), which
leads to labyrinthine patterns. Most schemes (RK4, ETDRK4,
ETDRK4-B) use a fixed time step of $\Delta t=0.1$, CK45 adapts its
time step to contain local absolute error below $10^{-4}$, leading
to an average $\Delta t \approx 0.62$. The smaller computational
time (about $11.6$ against $22.7$ minutes) is paid out with a
larger error. Note that Matlab is not very efficient when it comes
to loops; with Fortran the difference in execution times is wider
($2.5$ versus $8.7$ minutes). } \label{fig:error_comparison_laby}
\end{figure}

In figure (\ref{fig:error_comparison_laby}) a $2$-D comparison is
performed. As well as the previous schemes, we also employed the
Cash-Karp version of the RK4 scheme that adapts the timestep (the
overhead is small, easily allowing this). It is, again, clear that
the ETDRK4 and ETDRK4-B schemes are more accurate and over longer
times ETDRK4-B is the preferred scheme. The adaptive method is
very fast, and the accuracy can be improved by lessening the error
tolerances, and thus it is recommended for longer, more
time-consuming, computations.

\section{Concluding Remarks}

We have developed and packaged a suite of algorithms for solving
reaction diffusion equations. To make the algorithms immediately
relevant and directly usable for those in the reaction diffusion
equations community we have illustrated the algorithms upon recent
and varied examples from the literature. Probably the most
striking feature to emerge is how splendidly the method copes with
sharp variations in the solutions, and also how fast and accurate
the method is even with large timesteps.

But, before further congratulating ourselves upon the efficiency of spectral
methods we must discuss several disadvantages:

The scheme we present is utterly reliant upon the
reaction-diffusion equations being semi-linear, that is, the
diffusion terms are simply $u_{xx}+ C u_{yy}$ and similarly for
$v$; for some constant $C$ (we have simply had isotropic diffusion
in this article).

We have not discussed possible problems with aliasing, earlier
versions of our code utilized Orzag's 2/3 rule to filter this out.
However, this actually made no discernable difference to the
solutions and we later just discarded this. It is evident that the
method has terms $\exp(-\Omega^2\Delta t)$ so higher order modes
are, in any case, exponentially decaying; aliasing transfers some
lower order modes to higher ones, so for diffusion-like problems
the aliasing is automatically damped. Nonetheless aliasing is an
issue that should be borne in mind in any spectral scheme.

Fourier spectral methods require periodicity, and we are not in
the position, at least here, to set Neumann or Dirichlet boundary
conditions on the edge of the domain. That requires an extension
to Chebyshev, or some other basis functions. Thus we have to take
the domain size large enough that the waves, pulses, structures of
interest do not interact with the edges of the domain. In fact,
one can set up and indeed solve Dirichlet/ Neumann boundary
condition problems using integrating factor methods, see for
instance \cite{kassam03a}, but there is then an essential
difference. One must take the exponential of a full matrix, the
periodic case treated here is special as those matrices are then
diagonal and this simplification underlies all that we have done
here, and computing the exponential of a matrix is numerically
expensive particularly if it must be re-computed. This is an area
that deserves further thought and work as the prospective pay-off
in generating explicit timestepping codes for stiff PDEs in high
spatial dimensions is considerable.

In some cases spectral accuracy means that we can use so few modes
that the graphical solutions look unnaturally poor. This is
despite the isolated values at the grid points being spectrally
accurate, we can then utilize interpolation onto a finer grid
using periodic cardinal functions, an algorithm  for 1D is
supplied in \cite{weideman01a}. We present an alternative, and
generalization to 2D, in Appendix~\ref{app:cardinal} based upon
padding a Fourier transform with zeros.

Note that we are not claiming that the codes herein are the
absolute best algorithms available for reaction-diffusion
equations, nor do we attempt to imply that other scientists using
alternative algorithms have been misguided. In particular, in 1D,
the adaptive scheme of \cite{blom94a} has proved itself to be a
useful and accurate algorithm that we have enjoyed working with.
One alternative scheme that certainly suggests itself is an
Alternating Direction Implicit scheme where spatial discretization
is again done through spectral methods, for completeness we
provide a Matlab code that does this and we discuss this further
in an appendix. In essence, we find that the low-order time solver
usually used means that the scheme performs much less well than
the integrating factor method of the main text.

Our aim has been, and is, to provide good, clear, working,
versatile spectral schemes, that avoid stiffness issues, in a form
whereby they can be utilized and built upon by other scientists.
Thus, we hope, allowing them to concentrate upon the physics,
biology, chemistry or other scientific issue rather than upon
numerical concerns; the codes are summarized in
table~\ref{tab:table}, and are documented both internally and via
an electronic README file.

\begin{acmtable}{\textwidth}
\centering
\begin{verbatim}

Fortran                              Matlab
|-- OneD                             |-- OneD
|   |-- CK45                         |   |-- CK45
|   |   |-- auto_CK45.f              |   |   |-- auto_CK45.m
|   |   |-- epidemic_CK45.f          |   |   |-- epidemic_CK45.m
|   |   |-- fisher_CK45.f            |   |   |-- fisher1D_CK45.m
|   |   `-- gray1D_CK45.f            |   |   `-- gray1D_CK45.m
|   |-- ETDRK4_B                     |   |-- ETDRK4
|   |   |-- auto_ETDRK4_B.f          |   |   `-- gray1D_ETDRK4.m
|   |   |-- epidemic_ETDRK4_B.f      |   |-- ETDRK4_B
|   |   |-- fisher_ETDRK4_B.f        |   |   |-- auto_ETDRK4_B.m
|   |   `-- gray1D_ETDRK4_B.f        |   |   |-- epidemic_ETDRK4_B.m
|   `-- RK4                          |   |   |-- fisher1D_ETDRK4_B.m
|       |-- auto_RK4.f               |   |   `-- gray1D_ETDRK4_B.m
|       |-- epidemic_RK4.f           |   `-- RK4
|       |-- fisher_RK4.f             |       |-- auto_RK4.m
|       `-- gray1D_RK4.f             |       |-- epidemic_RK4.m
`-- TwoD                             |       |-- fisher1D_RK4.m
    |-- CK45                         |       `-- gray1D_RK4.m
    |   |-- gray2D_CK45.f            `-- TwoD
    |   `-- labyrinthe2D_CK45.f          |-- CK45
    |-- ETDRK4_B                         |   |-- fisher2D_CK45.m
    |   |-- gray2D_ETDRK4_B.f            |   |-- gray2D_CK45.m
    |   `-- labyrinthe2D_ETDRK4_B.f      |   `-- labyrinthe2D_CK45.m
    `-- RK4                              |-- ETDRK4
        |-- gray2D_RK4.f                 |   `-- labyrinthe2D_ETDRK4.m
        `-- labyrinthe2D_RK4.f           |-- ETDRK4_B
                                         |   |-- fisher2D_ETDRK4_B.m
Useful                                   |   |-- gray2D_ETDRK4_B.m
|-- adifisher.m                          |   `-- labyrinthe2D_ETDRK4_B.m
|-- fourierupsample.m                    `-- RK4
|-- fourierupsample2D.m                      |-- fisher2D_RK4.m
|-- plot_fisher2D.m                          |-- gray2D_RK4.m
|-- plot_gray2D.m                            `-- labyrinthe2D_RK4.m
`-- plot_labyrinthe2D.m

\end{verbatim}
\begin{tabular}{c|cc}
   & File Name & ref. Equation \\
  \hline
  \multirow{4}{0.08\linewidth}{1-D} & auto & (\ref{eqn:auto}) \\
      & epidemic & (\ref{eqn:epidemic}) \\
      & fisher1D & (\ref{eqn:fisher1D}) \\
      & gray1D & (\ref{eqn:gray1D}) \\
  \hline
  \multirow{3}{0.08\linewidth}{2-D} & fisher2D & Appendix~\ref{app:adi}\\
      & gray2D & (\ref{eqn:gray2D}) \\
      & labyrinthe2D & (\ref{eqn:labyrinthe2D}) \\
  \hline
\end{tabular}
\vspace{0.5cm}
\caption{A schematic tree of the provided algorithms. Names and
corresponding equations are matched in the lower table.}
\label{tab:table}
\end{acmtable}

\clearpage
\appendix
\section*{APPENDIX}

\section{The one dimensional Matlab code}
\label{app:code}

\begin{verbatim}
function gray1D_RK4(N,Nfinal,dt,ckeep,L,epsilon,a,b)
if nargin<8;
    disp('Using default parameters');
    N=512; Nfinal=10000; dt=0.2; ckeep=10;
    L=50; epsilon=0.01; a=9*epsilon; b=0.4*epsilon^(1/3);
end

x=(2*L/N)*(-N/2:N/2-1)';
u=initial(x,L); uhat=fft(u);
ukeep=zeros(N,2,1+Nfinal/ckeep);
ukeep(:,:,1)=u;
tkeep=dt*[0:ckeep:Nfinal];
ksq=((pi/L)*[0:N/2 -N/2+1:-1]').^2;
%-----------------Runge-Kutta----------------------------------
E=[exp(-dt*ksq/2) exp(-epsilon*dt*ksq/2)]; E2=E.^2;
for n=1:Nfinal
    k1=dt*fft(rhside(u,a,b));
    u2=real(ifft(E.*(uhat+k1/2)));
    k2=dt*fft(rhside(u2,a,b));
    u3=real(ifft(E.*uhat+k2/2));
    k3=dt*fft(rhside(u3,a,b));
    u4=real(ifft(E2.*uhat+E.*k3));
    k4=dt*fft(rhside(u4,a,b));
    uhat=E2.*uhat+(E2.*k1+2*E.*(k2+k3)+k4)/6;
    u=real(ifft(uhat));
    if mod(n,ckeep)==0,
        ukeep(:,:,1+n/ckeep)=u;
    end
end
save('gray1D_RK4.mat','tkeep','ukeep','N','L','x')
%----------------------Figures---------------------------------
mesh(tkeep,x,squeeze(ukeep(:,2,:))); view([60,75]);
xlabel('t'); ylabel('x'); zlabel('z');
title('(a) Surface plot of v')
%--------------Initial Condition ------------------------------
function u=initial(x,L)
u=[1-0.5*(sin(pi*(x-L)/(2*L)).^100) ...
    0.25*(sin(pi*(x-L)/(2*L)).^100)];
%---------------Right Hand Side--------------------------------
function rhs2=rhside(u,a,b)
t1=u(:,1).*u(:,2).*u(:,2);
rhs2=[-t1+a*(1-u(:,1)) t1-b*u(:,2)];
\end{verbatim}
This produces figure \ref{fig:gscott1D}(a) of the text for the
Gray-Scott equations.

\section{Fourier cardinal function interpolation}
\label{app:cardinal}

As noted in the text, spectral methods are often very accurate
even with few interpolation points. When plotting graphically this
sometimes leads to artificially poor-looking output, clearly the
solution is spectrally accurate at each interpolation point and we
just need to insert more points. Fourier cardinal function
interpolation as in \cite{weideman01a} can be used in 1D, or,
more in tune with the current article, one can pad an FFT with
additional zeros and then invert which is convenient in either one
or two space dimensions. The short Matlab scripts that do this
are:

In one dimension:
\begin{verbatim}
function fout=fourierupsample(fin,newN);
% Given a periodic function fin, computed at N equispaced nodes in
% the periodic domain [-L,L], fout is its upsampled version on newN
% nodes onto the same domain.

N=length(fin); HiF=(N-mod(N,2))/2+1;
fftfin=max((newN/N),1)*fft(fin);
fout=real(ifft([fftfin(1:HiF); zeros(newN-N,1); fftfin(HiF+1:N)]));
\end{verbatim}

And in two dimensions:

\begin{verbatim}
function fout=fourierupsample2D(fin,newNx,newNy);
% Given a periodic function fin in 2D computed at equidistant
% nodes Nx x Ny, then fout is its upsampled version on
% newNx x newNy nodes.

[Ny,Nx]=size(fin);
HiFx=(Nx-mod(Nx,2))/2+1; HiFy=(Ny-mod(Ny,2))/2+1;
fftfin=max((newNx/Nx)*(newNy/Ny),1)*fft2(fin);
fout=real(ifft2([fftfin(1:HiFy,1:HiFx), ...
    zeros(HiFy,newNx-Nx), fftfin(1:HiFy,HiFx+1:end); ...
    zeros(newNy-Ny,newNx); fftfin(HiFy+1:end,1:HiFx), ...
    zeros(Ny-HiFy,newNx-Nx), fftfin(HiFy+1:end,HiFx+1:end)]));
\end{verbatim}

\section{Alternating Direction Implicit (ADI) Methods: }
\label{app:adi}

\begin{figure}
\leavevmode
\centering
\includegraphics[height=5cm]{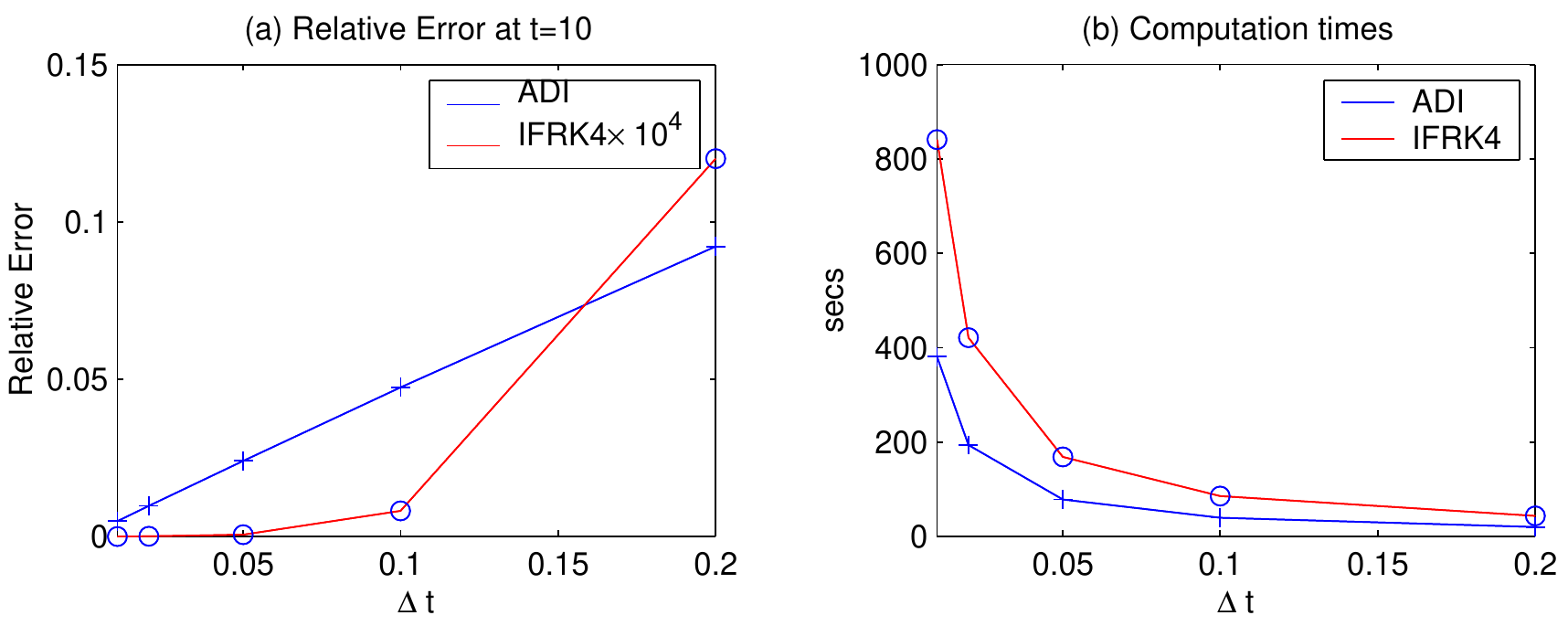}
\vspace*{2mm}
\caption{Panels (a) and (b) relative errors and timings for ADI
versus the integrating factor method.} \label{fig:adi}
\end{figure}

This appears to be a viable alternative to that which we have
presented in the main text; it is worth briefly outlining the
method. The basic idea is similar to operator (Strang) splitting
and the method is discussing in some detail in
\cite{boyd01a,press92a}. As noted by Boyd the conventional
centered finite difference schemes can easily be modified by using
spectral differentiation matrices, let $D$ denote the $N\times N$
Fourier differentiation matrix (\cite{fornberg98a},
\cite{boyd01a},\cite{trefethen01a}\cite{weideman01a}). ADI, or at
least the ADI we use here, means that we split each time step into
two and we first deal implicitly with one set of space derivatives
and then in the next half time step with the other. So
\[
  u_t=u_{xx}+u_{yy}+f(u)
\]
is approximated by the following matrix system, here $I$ is the
$N\times N$ identity matrix and $U$ is the matrix of $u$
\begin{eqnarray*}
U^{(t+1/2)}&=&\left[I-\frac{\Delta t}{2}D\right]^{-1}
U^{(t)}\left[I+\frac{\Delta t}{2}D^T\right]+
\left[I-\frac{\Delta t}{2}D\right]^{-1}\frac{\Delta t}{2} F(u^{(t)})
\\
U^{(t+1)}&=&\left[I+\frac{\Delta t}{2}D\right]
U^{(t+1/2)}\left[I-\frac{\Delta t}{2}D^T\right]^{-1}+
\frac{\Delta t}{2} F(u^{(t+1/2)})\left[I-\frac{\Delta t}{2}D^T\right]^{-1}
\end{eqnarray*}
with the evident advantage that each matrix is evaluated only once
and thereafter we just have matrix multiplication. This is ideal
for implementation in Matlab. Unfortunately, this is only accurate
to $O(\Delta t)^2$ so although the method is nicely stable one
requires relatively small time steps relative to the explicit
non-stiff scheme that is used in the main text. For instance, for
Fisher's equation in 2D we show some comparative errors and
timings in figure \ref{fig:adi}; note this simulation is for
Fisher's equation using $256\times 256$ Fourier modes on a
$50\times 50$ domain, the initial condition is a Gaussian
$0.2\exp(-0.25(x^2+y^2))$. The integrating factor solution with
$\Delta t=10^{-2}$ is taken as the reference solution and the
relative errors are computed as the maximal difference away from
this. Notably the errors from the integrating factor scheme are
multiplied by $10^4$ in order that they are visible.

Doubtless one could improve the naive implementation above, but
for the application to reaction-diffusion equations it seems
uncompetitive. There are advantages though, in that it is
generalizable to problems with $u$ dependent diffusivity whereas
the integrating factor method is not.

\bibliography{fronts}
\bibliographystyle{acmtrans}

\end{document}